\documentclass[reqno,12pt]{amsart}
\usepackage[foot]{amsaddr}
\usepackage{amssymb}
\usepackage{caption}
\usepackage{graphicx}
\usepackage{subcaption}
\usepackage{threeparttable}
\usepackage{mathrsfs}
\usepackage{float}
\usepackage{dsfont}
\usepackage{amsmath}
\usepackage{booktabs}
\usepackage{array}
\usepackage{multirow}
\usepackage{lineno}
\usepackage{epsfig}
\usepackage{makecell}
\usepackage{diagbox}
\usepackage{changepage}
\usepackage[figuresright]{rotating}
\usepackage[colorlinks=true]{hyperref}
\usepackage{setspace}
\hypersetup{urlcolor=red, citecolor=blue}
\usepackage{cite}
\usepackage{enumerate}
\usepackage{appendix}

\usepackage[dvipsnames]{xcolor}
\usepackage{tikz}

\allowdisplaybreaks[2]
\newtheorem{theorem}{Theorem}[section]

\newtheorem{lemma}[theorem]{Lemma}
\newtheorem{proposition}[theorem]{Proposition}
\newtheorem{definition}[theorem]{Definition}
\newtheorem{remark}{Remark}[section]

\newtheorem{example}{Example}[section]
  
  \newtheorem{assumption}{Assumption}[section]
  \DeclareMathOperator{\dive}{div}
\numberwithin{equation}{section}
\title[Non-convergence of the principal eigenvalue]{Non-convergence of the principal eigenvalue of elliptic operators for large advection}

\begin{document}

\begin{abstract}
This paper investigates the limit of the principal eigenvalue $\lambda(s)$ as $s\to+\infty$ for the following elliptic equation
 \begin{align*}
-\Delta\varphi(x)-2s\mathbf{v}\cdot\nabla\varphi(x)+c(x)\varphi(x)=\lambda(s)\varphi(x), \quad x\in \Omega
 \end{align*}
in a bounded domain $\Omega\subset \mathbb{R}^d (d\geq 1)$  with the Neumann boundary condition.
Previous studies have shown that under certain conditions on $\mathbf{v}$, $\lambda(s)$ converges as $s\to\infty$ (including cases where $\lim\limits_{s \to\infty }\lambda(s)=\pm\infty$). This work constructs an example such that $\lambda(s)$ is divergent as $s\to+\infty$. This seems to be the first rigorous result demonstrating the non-convergence of the principal eigenvalue for second-order linear elliptic operators with some strong advection. As an application, we demonstrate that for the classical advection-reaction-diffusion model with advective velocity field $\mathbf{v}=\nabla m$, where $m$ is a potential function with infinite oscillations, the principal eigenvalue  changes sign infinitely often along a subsequence of $s\to\infty$. This leads to   solution behaviors that differ significantly from those observed when $m$ is non-oscillatory.

\end{abstract}

\renewcommand{\thefootnote}{\fnsymbol{footnote}}
\author[X, Bai, Z.A. Wang, X. Xu, K. Zhang, M. Zhou]{Xueli Bai$^{*}$, Zhi-An Wang$^{\dagger}$, Xin Xu$^{\dagger}$, Kexing Zhang$^{\S}$, Maolin Zhou$^{\sharp}$ }
\footnotetext{$^{*}$ School of Mathematics and Statistics, Northwestern Polytechnic University, Xi'an 710072, China;  xlbai2015@nwpu.edu.cn}
\footnotetext{$^{\dagger}$ Department of Applied Mathematics, The Hong Kong Polytechnic University, Hung Hom, Hong Kong, P.R. China; mawza@polyu.edu.hk}
\footnotetext{$^{\ddagger}$ School of Mathematical Sciences, South China Normal University, Guangdong 510631, China; xuxin1994pkq@163.com}
\footnotetext{$^{\S}$ Chern Institute of Mathematics and LPMC, Nankai University, Tianjin 300071, China; kxzmath@163.com}
\footnotetext{$^{\sharp}$ Chern Institute of Mathematics and LPMC, Nankai University, Tianjin 300071, China; zhouml123@nankai.edu.cn}

\maketitle

\section{Introduction}
In this paper, we investigate the limit of the principal eigenvalue $\lambda(s)$, as $s\to+\infty$, for the following second-order elliptic Neumann eigenvalue problem:
\begin{equation}\label{CF1}
-\Delta\varphi(x)-s\mathbf{v}\cdot\nabla\varphi(x)+c(x)\varphi(x)=\lambda(s)\varphi(x),\ x\in \Omega,
\end{equation}
where $\Omega\subset\mathbb{R}^d\ (d\geq 1)$ is a bounded domain with smooth boundary,  $\mathbf{v}\in C(\bar{\Omega})$ is a drift (or advective) velocity field, $c\in C(\bar{\Omega})$ is a given function.

The asymptotic behavior of the principal eigenvalue and its eigenfunctions for second-order linear elliptic operators has been an important research topic \cite{CCL6,CCL7}. When the Dirichlet boundary condition is prescribed to \eqref{CF1} with $c(x)\equiv 0$, the asymptotic behaviors of the principal eigenvalue $\lambda(s)$ as $s \to +\infty$ (or equivalently, the coefficient of $\Delta \varphi$ tends to zero) were extensively discussed in \cite{DEF73, FA7273, WENTZEL, W75} under various conditions in the vector field $\mathbf{v}$, followed by work \cite{DF78, FS97} on the asymptotic behavior of the principal eigenfunctions.

When the vector field $\mathbf{v}$ is divergence-free (i.e., $\dive{\mathbf{v}}=0$), Berestycki et al. in \cite{Bcmp} introduced the concept of the first integral of $\mathbf{v}$ ($w$ is said to be a first integral of $\mathbf{v}$ if $w\in H^1(\Omega)$, $w\neq 0$, and $\mathbf{v}\cdot\nabla w=0 $ almost everywhere in $\Omega$, see \cite[Definition 0.1]{Bcmp}). They showed that the principal eigenvalues $\lambda(s)$, subject to Dirichlet boundary conditions, are bounded if and only if $\mathbf{v}$ has the first integral in $H_0^1(\Omega)$, where $\lambda(s)$ converges to the minimum value of Rayleigh quotient of a first integral, if any. The above striking results also hold for Neumann or periodic boundary conditions with additional condition $\textbf{v}\cdot \textbf{n}=0$ on $\partial \Omega$. The only difference is that the set of first integrals of $\mathbf{v}$ is not empty for Neumann or periodic boundary conditions, since it contains all nonzero constant functions.
Later, Chen and Lou in \cite{CL} initiated research for the vector field that is not divergence-free by assuming that $\textbf{v}$ is a gradient field: $\mathbf{v}=\nabla m$, where $m(x)\not \equiv 0$ is a scalar potential function in $\Omega$. Specifically, they considered the following Neumann eigenvalue problem
\begin{equation}\label{CF}
  \left\{
\begin{array}{llll}
  -\Delta\varphi(x)-2s\nabla m\cdot\nabla\varphi(x)+c(x)\varphi(x)=\lambda(s)\varphi(x),& x\in \Omega,\\[2mm]
  \nabla\varphi\cdot\mathbf{n}=0,&x\in\partial\Omega,
  \end{array}\right.
\end{equation}
where $\mathbf{n}=\mathbf{n}(x)$ denotes the unit exterior normal vector at $x\in \partial\Omega$,
and proved that if the critical points of $m$ are nondegenerate, the principal eigenvalue $\lambda(s)$ converges as $s \to +\infty$ and satisfies
  \begin{align}\label{clr}
  \lim_{s\to+\infty}\lambda(s)=\min_{\mathcal{M}}c(x),
  \end{align}
where $\mathcal{M}$ denotes the set of points of local maximum of $m(x)$.
In one-dimensional space, if $m(x)$ vanishes on a sub-interval $[a,b]\subset[0,1]$ with $0<a<b<1$,  and the derivative $m^\prime(x)$ changes sign finitely many times, Peng and Zhou in \cite{PZ} obtained the convergence of $\lambda(s)$.
In \cite{LLPZ}, Liu et al. demonstrated that for a general vector field $\mathbf{v}$ in two dimensions, the principal eigenvalue $\lambda(s)$ converges as $s\to+\infty$, and the asymptotic behaviors are completely determined by some connected components in the omega-limit set of the system of ordinary differential equations associated with the drift term, which includes stable fixed points, stable limit cycles, hyperbolic saddles connecting homoclinic orbits, and families of closed orbits.
Recently, the asymptotics of the principal eigenvalues of second-order parabolic operators with small diffusion or large advection were investigated in the literature (cf. \cite{HNR2011,LLPZd,LLPZx,PZhao,N09,N10} and references therein).

The existing results, among others, aim mainly to find conditions in the vector field so that the principal eigenvalue $\lambda(s)$ converges as $s\to+\infty$, as mentioned above. However, it remains elusive whether there is a vector field $\mathbf{v}$ such that the principal eigenvalue $\lambda(s)$ is divergent as $s\to\infty$. The goal of this paper is to attempt this question and give a positive answer to the Neumann eigenvalue problem \eqref{CF}.
\par
\vspace{0.2cm}

Our main result is stated in the following theorem.
\begin{theorem}\label{robin}
Let $\Omega=B_1(0)\subset \mathbb{R}^d$. Then there exist potential functions $m\in C^1(\Omega)$ that vanish on some subdomains of $\Omega$ such that the corresponding principal eigenvalues $\lambda(s)$ satisfying \eqref{CF} diverge as $s\to+\infty$. Specifically, for each principal eigenvalue, we can find two convergent subsequences $\{s_i\}_{i=1}^{+\infty}$ and $\{s_j\}_{j=1}^{+\infty}$ with $s_i\to+\infty$  and $s_j\to+\infty$, such that $\lim\limits_{s_i\to+\infty}\lambda(s_i)$ and $\lim\limits_{s_j\to+\infty}\lambda(s_j)$ are two different constants.
\end{theorem}

\begin{remark}\em{
As far as we know, the result in Theorem \ref{robin} seems to be the first one rigorously showing the non-convergence of the limit of principal eigenvalue for second-order linear elliptic operators for large advection. The potential functions constructed in Theorem \ref{robin} are required to exhibit infinite oscillations about the $r$-axis with non-periodic recurrence near the boundary of the degenerate subdomain (i.e., the subdomain where $m$ vanishes). The amplitude of these oscillations decays as they approach the degenerate subdomain. The detailed construction of such $m$ will be provided in Sections 2 and 5. }
\end{remark}

\begin{remark}\em{
As an application of Theorem \ref{robin}, we can show that the nonlinear reaction-diffusion-advection equation, where the advection vector $\mathbf{v}$ is a gradient field (i.e., $\mathbf{v}=\nabla m$ for some potential function $m$), may demonstrate an interesting phenomenon: there exist two distinct subsequences of $s\to+\infty$, along which the system alternates ecological regimes of persistence and extinction (see Section \ref{S6}).
}
\end{remark}

The rest of this paper is organized as follows. In Section \ref{prel}, we provide some preliminary results and introduce some notions. In Sections \ref{dddnnn}, we characterize the limit of the principal eigenvalue for the gradient vector field $\mathbf{v}=\nabla m(x)$ with two different types of $m$, respectively.
Then we present the proof of Theorem \ref{robin} in Section \ref{notconver}, utilizing the results prepared in the preceding sections. In Section \ref{S6}, we show that for reaction-diffusion-advection models with the advection terms constructed in Section \ref{notconver}, the principal eigenvalue can change signs infinitely often along two subsequences of $s\to\infty$, as an application of Theorem \ref{robin}.

\section{Preliminaries }\label{prel}
In this section, we introduce some notions and properties of relevant eigenvalue problems, focusing on the radial symmetric domain $\Omega=B_{1}(0)\subset \mathbb{R}^d$ and the case where $\mathbf{v}$ is a gradient vector field, i.e.,  $\mathbf{v}(x)=\nabla m(x)$.

\subsection{Variational characterization and properties of principal eigenvalues}

Assume that $m(x)$ and $c(x)$ are radially symmetric functions. We write $m(r)=m(|x|)$ and $c(r)=c(|x|)$, where $r=|x|$.
Then, equation \eqref{CF} becomes
 \begin{equation}
\left\{
\begin{array}{l}
-\varphi''(r)-\frac{d-1}{r}\varphi'(r)-2sm'(r)\varphi'(r)+c(r)\varphi(r)=\lambda(s)\varphi(r),\quad 0<r<1,\\
\varphi'(0)=\varphi'(1)=0,
\end{array}\right.\label{eq}
\end{equation}
where $\lambda(s)$ is the principal eigenvalue.
Let us begin by recalling the variational characterization of the principal eigenvalue of problem \eqref{eq}.

{\bf Variational characterization of $\lambda(s)$.} It is well known that the principal eigenvalue $\lambda(s)$ can be characterized by ($\!\!$\cite{Evans}):
\begin{eqnarray}
\lambda(s)
&=&\min_{\varphi\in H^1((0,1))\atop\int_0^1 r^{d-1}e^{2sm(r)}\varphi^2\mathrm{d}r=1  }\int_0^1 r^{d-1}e^{2sm(r)}\left(|\varphi'|^2+c(r)|\varphi|^2\right)\mathrm{d} r\label{voe}\nonumber\\
&=&\min_{w\in H^1((0,1))\atop\int_0^1 w^2\mathrm{d} r=1 }\int_0^1 \left|w'-sm'(r)w-\frac{d-1}{2r}w\right|^2+c(r)|w|^2\mathrm{d} r.\label{voa}
\end{eqnarray}
Moreover, the minimizers of \eqref{voe} are the principal eigenfunction of \eqref{eq}, denoted by $\varphi_1(s;r)$ and $w_1(s;r)$, where $w_1$ satisfies the normalization condition
$\int_0^1 w_1^2\mathrm{d} r=1$
and is related to $\varphi_1(s;r)$ through the transformation $w_1(s;r)=e^{sm(r)}r^{\frac{d-1}{2}}\varphi_1(s;r).$

Let $c_*=\min_{r\in[0,1]}c(r)$ and $c^*=\max_{r\in[0,1]}c(r)$. In the following, we assume $c_*> 0$. The case when $c_*\leq 0$ can be treated similarly, as the proof essentially follows the same argument with only minor modifications. From the variational formula \eqref{voe}, one can directly see that $\lambda(s)$ is bounded below by $c_*$ and above by $c^*$, namely
\begin{equation}
c_*\le \lambda(s)\le c^*,\quad \forall\ s>0.\label{ubd}
\end{equation}

Hereafter, we will always assume that $m\equiv 0$ for $ r\in[a,b]$ with $0<a<b<1$.  The following equation
\begin{align}\label{AB}
-\varphi''(r)-\frac{d-1}{r}\varphi'(r)+c(r)\varphi(r)=\lambda\varphi(r),\quad r\in\left(a,b\right).
\end{align}
will be frequently used.  For convenience, we introduce the following definition:
\begin{definition}\label{d2}
Let $\lambda^{\mathcal{D}}$ be the principal eigenvalue of \eqref{AB} subject to the Dirichlet boundary condition: $\varphi(a)=\varphi(b)=0$. We denote by $\varphi^{\mathcal{D}}(r)$ for $r\in (a,b)$ the principal eigenfunction corresponding to the principal eigenvalue $\lambda^{\mathcal{D}}$ with normalization $$\int_a^br^{d-1}|\varphi^{\mathcal{D}}(r)|^2\mathrm{d} r=1.$$
Similarly, we use  $\lambda^{\mathcal{N}}$ and $\varphi^{\mathcal{N}}(r)$ to denote the principal eigenvalue and eigenfunction of \eqref{AB} with Neumann boundary condition: $\varphi'(a)=\varphi'(b)=0$.
\end{definition}

{\bf Variational characterization of $\lambda^{\mathcal{D}}$ and  $\lambda^{\mathcal{N}}$.} The variational characterizations of $\lambda^{\mathcal{D}}$ and  $\lambda^{\mathcal{N}}$ are given by
\begin{equation}\label{vd}
\lambda^{\mathcal{D}}=\min_{\varphi\in H_0^1((a,b))\atop\int_a^br^{d-1}\varphi^2\mathrm{d} r=1}\int_a^br^{d-1}(|\varphi'|^2+c(r)|\varphi|^2)\mathrm{d} r,
\end{equation}
and
\begin{equation}\label{vn}
\lambda^{\mathcal{N}}=\min_{\varphi\in H^1((a,b))\atop\int_a^br^{d-1}\varphi^2\mathrm{d} r=1}\int_a^br^{d-1}(|\varphi'|^2+c(r)|\varphi|^2)\mathrm{d} r,
\end{equation}
where we refer to \cite{Evans}.
Moreover, by the maximum principle \cite{GT}, we have
\begin{equation}\label{DDN}
\lambda^{\mathcal{D}}>\lambda^{\mathcal{N}}.
\end{equation}

\subsection{Assumptions and notions}

In this subsection, we will define two different types of step functions to aid in describing the construction of $m(r)$ on $[\delta,a]\cup[b,1-\delta]$, where $\delta\in(0,a)$ is a fixed constant. Before proceeding with the definitions, we impose the following  assumptions on the non-constant functions $m(\cdot)\in C^1([0,1])$ and $c(\cdot)\in C([0,1])$:\par
\begin{itemize}
\item[(H1)] $m(r)=m(1-r)$ on $[0,1]$, and $m(r)=0$ on $[a,b]$ with $0<a<b<1$, where $a+b=1$.\par
\item[(H2)] $c(r)>0$ on $[0,1]$ (which implies $\lambda^{\mathcal{D}}>0$) and $c(x)>\lambda^{\mathcal{D}}$ for $r\in [0,a]\cup[b,1].$ \par
\end{itemize}

\begin{remark}\label{*}
Hypothesis {\rm{(H1)}} is employed to simplify the computation technicalities. For the non-symmetric function $m(r)$, we can apply a calculation method same as that used for the symmetric case. By selecting $c(r)$ to be suitably small on the interval $(a,b)$, we can identify a function $c(r)$ that satisfies $c(r)>\lambda^{\mathcal{D}}$ for $r\in[0,a]\cup[b,1]$.  Therefore, the hypothesis  {\rm(H2)} is valid.

\end{remark}

Now, we can define two different types of step functions.

\begin{definition}\label{D3}
Let $\delta\in(0,a)$ and $h,\ \alpha,\ \beta,\ \nu,\ l $  be constants satisfying  $0<h<\alpha<\beta<1<\nu$, $l\in\mathbb{N}$ and $\sum_{i=1}^{+\infty}(\alpha^{i+l}+\beta^{i+l})=a-\delta$.

{\rm(1)} We define $\tilde{m}(\delta,h,\alpha,\beta,\nu,l;r)$ for $r\in [\delta,a]$ as follows:
\begin{equation}\label{tilm}
\tilde{m}(\delta,h,\alpha,\beta,\nu,l;r) = \left\{ \begin{array}{lll}
-h^{n}, &r\in [x_n, y_n),&n\in \mathbb{N},\\
\nu h^{n}, &r\in [ y_n,x_{n+1}),&n\in \mathbb{N},
 \end{array} \right.
\end{equation}
where $x_0=\delta$, $y_0=\delta+\alpha^{1+l}$, $x_n=\delta+\sum_{i=1}^{n}(\alpha^{i+l}+\beta^{i+l})$, and $y_n=\delta+\sum_{i=1}^{n}(\alpha^{i+l}+\beta^{i+l})+\alpha^{n+l+1}$ for $n\ge1$.
For $r\in[b,1-\delta]$, we define $\tilde{m}(\delta,h,\alpha,\beta,\nu,l;r)=\tilde{m}(\delta,h,\alpha,\beta,\nu,l;1-r)$.
\par
{\rm(2)}
We define $\bar{m}(\delta,h,\alpha,\beta,\nu,l;r)$ for $r\in [\delta,a]$ as follows:
\begin{equation}\label{barm}
\bar{m}(\delta,h,\alpha,\beta,\nu,l;r) = \left\{ \begin{array}{lll}
h^{n},& r\in[Y_{n-1},X_n],&n\in\mathbb{N},\\
-\nu h^{n},&r\in[X_n,Y_n], &n\in\mathbb{N},
 \end{array} \right.
\end{equation}
where
$Y_0=\delta$, and
$X_n=\delta+\sum_{i=1}^n (\alpha^{i+l}+\beta^{i+l})-\beta^{n+l}$,
$Y_n=\delta+\sum_{i=1}^n (\alpha^{i+l}+\beta^{i+l})$ for $n\ge1$.
For $r\in[b,1-\delta]$, we define $\bar{m}(\delta,h,\alpha,\beta,\nu,l;r)=\bar{m}(\delta,h,\alpha,\beta,\nu,l;1-r)$.
\end{definition}

Based on the definitions of $\tilde{m}(\delta,h,\alpha,\beta,\nu,l;r)$ and $\bar{m}(\delta,h,\alpha,\beta,\nu,l;r)$, we construct the following two sets.

\begin{definition}\label{e1}
Let $m(r)\in C^1([0,1])$ and $m'(r)$ changes sign finitely many times in $[0,\delta)\cup(1-\delta,1]$.
\begin{itemize}
\item[(1)] {\rm(The set $S_{\mathcal{D}}$\rm)} We say $m(r)\in S_{\mathcal{D}}$ if there exist constants $\delta$, $h$, $\alpha$, $\beta$, $\nu$, $l$ satisfying $0<h<\alpha< \beta<1<\nu$ and $l\in\mathbb{N}$ such that for any $r\in [\delta,a]\cup[b,1-\delta]$, $m(r)\ge \tilde{m}(\delta,h,\alpha,\beta,\nu,l;r)$.\par
\item[(2)] {\rm(The set $S_{\mathcal{N}}$\rm)} We say
$m(r)\in S_{\mathcal{N}}$ if there exist constants $\delta$, $h$, $\alpha$, $\beta$, $\nu$, $l$ satisfying $0<h<\alpha<\beta<1<\nu$ and $l\in\mathbb{N}$ such that for $r\in[\delta,a]\cup[b,1-\delta]$, $m(r)\le \bar{m}(\delta,h,\alpha,\beta,\nu,l;r)$.
\end{itemize}
\end{definition}
To understand Definition \ref{e1}, we provide the following example, which will be applied to prove Theorem \ref{robin} in Section 5.
\begin{example}\label{EE}
Choose $\delta=x_0=\frac{1}{84}$, $h=\frac{1}{10}$, $\alpha=\frac{1}{8}$, $\beta=\frac{1}{4}$, $\nu=2$, $l=0$, as well as $a=\frac{41}{84}$, $b=\frac{43}{84}$. For $n\in\mathbb{N}$, we define
\begin{equation*}
m_{\mathcal{D}}(r)=\left\{
\begin{array}{lll}
20,& r\in[0,x_0),\\
h^{n-1}(\cos{\pi\frac{r-x_n}{z_n-x_n}}+1),& r\in[x_{n},z_n),\\
h^{n}(\cos{\pi(\frac{r-z_n}{y_n-z_n}-1)}+1),& r\in[z_{n},y_n),\\
2h^{n},& r\in[y_{n},x_{n+1}),\\
0, & r\in[a,b],\\
m_D(1-r), & r\in(b,1],
\end{array}\right.
\end{equation*}
where $x_n=\frac{41}{84}-\frac{1}{8^n}(\frac{1}{7}+\frac{2^n}{3})$, $y_n=\frac{41}{84}-\frac{1}{8^n}(\frac{1}{56}+\frac{2^n}{3})$ and $z_n=\frac{x_n+y_n}{2}$.
Let $\bar{m}$ be as defined in Definition \ref{D3} (1). It is routine to check that $m_{\mathcal{D}}(r)\ge \tilde{m}(\frac{1}{84},\frac{1}{10},\frac{1}{8},\frac{1}{4},2,0;r), r\in [\delta,a]\cup[b,1-\delta]$, which implies that $m_{\mathcal{D}}(r)\in S_{\mathcal{D}}$. Moreover, the function $m_{\mathcal{D}}$ satisfies properties $m_{\mathcal{D}}(r)> 0$ on $(x_n,z_n)\cup(z_n,y_n)$ and  $m_{\mathcal{D}}(z_n)=({m_{\mathcal{D}}})'(z_n)=0$.

 Relying on the definition of $m_{\mathcal{D}}$, we can define
\begin{equation*}
m_{\mathcal{N}}(r)=\left\{
\begin{array}{llll}
&m_{\mathcal{D}}(r),& r\in[0,z_{n_0}),\\
&-m_{\mathcal{D}}(r),& r\in[z_{n_0},a),
\end{array}\right.
\end{equation*}
where $n_0\in \mathbb{N}$. Take $\delta=Y_0=x_{n_0+1}$,  $h=\frac{1}{10^{n_0+1}}$, $\alpha=\frac{1}{8}$, $\beta=\frac{1}{4}$, $\nu=2$, $l=n_0+1$.
It is easy to see that $m_{\mathcal{N}}(r)\le \bar{m}(x_{n_0+1},\frac{1}{10^{n_0+1}},\frac{1}{8},\frac{1}{4},2,n_0+1;r)$ on $[\delta,a]\cup[b,1-\delta]$, where $\bar{m}$ is defined as in Definition \ref{D3} (2). Thus, we conclude that $m_{\mathcal{N}}(r)\in S_{\mathcal{N}}$ immediately.
\end{example}

\section{Convergence of $\lambda(s)$ as $s\to+\infty$ for $m(r)\in S_{\mathcal{D}}$ and $S_{\mathcal{N}}$}\label{dddnnn}

We  begin by demonstrating  the convergence property of $\varphi_1(s;r)$ on the degenerate interval $(a,b)$ as $s\to+\infty$, which will be used later.


\begin{lemma}\label{converge}
Suppose that $m(r)=0$ on the interval $[a,b]$ and  $\lambda_*=\liminf\limits_{s\to+\infty}\lambda(s)$. Then,
there exists a sequence $\{s_i\}_{i=1}^{+\infty}$ with $s_i\to+\infty$ as $i\to+\infty$, and a function $\varphi_*(r)$ such that $\varphi_1(s_i;r)\rightharpoonup \varphi_*(r)$ weakly in $H^1((a,b))$ and $\varphi_1(s_i;r)\to \varphi_*(r)$ in $C([a,b])\cap C^1_{\mathrm{loc}}((a,b))$ as $i\to+\infty$. Moreover, $\varphi_*(r)$ satisfies
\begin{equation}
-\varphi_*''(r)-\frac{d-1}{r}\varphi_*'(r)+c(r)\varphi_*(r)=\lambda_*\varphi_*(r),\quad r\in (a,b),\label{v*}
\end{equation}
in the weak sense.
\end{lemma}

\begin{proof}
For simplicity, let $C$ represent some constant(possibly different) that is independent of the essential variables. We first select the sequence $\{s_i\}_{i=1}^{+\infty}$ such that $\lim_{i\to+\infty}\lambda(s_i)=\lambda_*$.
In light of $m(r)=0$ for $r\in[a,b]$, \eqref{ubd}, and the variational characterization \eqref{voe}, we have
\begin{align*}
c^*&\ge \lambda(s_i)\ge a^{d-1}\int_a^b  |\varphi_1'(s_i;r)|^2+c(r)|\varphi_1(s_i;r)|^2\mathrm{d} r
\ge C\int_a^b |\varphi_1'(s_i;r)|^2+|\varphi_1(s_i;r)|^2\mathrm{d} r
\end{align*}
with $C=a^{d-1}\min\{1,c_*\}>0$. Hence, $\varphi_1(s_i;r)$ is bounded in $H^1((a,b))$ and we can select a subsequence of $\{s_i\}$, still denoted by itself, such that $\varphi_1(s_i;r)\rightharpoonup \varphi_*(r)$ weakly in $H^1((a,b))$ for some $\varphi_*(x)$.
By the compact embedding theorem, after passing to a subsequence, there exists a subsequence of $\{s_{i}\}_{i=1}^{+\infty}$ such that
\begin{equation}
\varphi_1(s_{i};r)\to \varphi_*(r), \quad\text{ in }\ C([a,b]),\quad i\to+\infty,\label{m}
\end{equation}
and a constant $C>0$, independent of $i$, such that
\begin{align}\label{wb}
|\varphi_1(s_{i};r)|\le C,\quad r\in[a,b].
\end{align}
Moreover, from equation \eqref{eq}, we have
\begin{equation*}\label{d}
-\varphi_1''(s_{i};r)-\frac{d-1}{r}\varphi_1'(s_i;r)+c(r)\varphi_1(s_{i};r)=\lambda(s_{i})\varphi_1(s_{i};r),\quad r\in(a,b).
\end{equation*}
Therefore, by a standard compactness consideration,
we can further extract a subsequence of $\{s_{i}\}_{i=1}^{+\infty}$, still denoted by itself, such that $\varphi_1(s_{i};r)\to\varphi_*(r)$ in $C^1_{\mathrm{loc}}((a,b))$ and  $\varphi_*(r)$ satisfies
\begin{equation*}
-\varphi_*''(r)-\frac{d-1}{r}\varphi_*'(r)+c(r)\varphi_*(r)=\lambda_*\varphi_*(r),
\end{equation*}
in the weak sense. The proof is complete.
\end{proof}

According to Definition \ref{D3},  
for simplicity, we only consider the case $l=0$ and the case $l\in\mathbb{N}^+$ can be proved in a similar manner. 

\subsection{Consider $m(r)\in S_{\mathcal{D}}$. }\label{DD}
In this section, we present the convergence property for $m(r)\in S_{\mathcal{D}}$.

\begin{proposition}\label{DT}
 If $m(r)\in S_{\mathcal{D}}$, then the principal eigenvalue $\lambda(s)$ of equation \eqref{eq} satisfies $\lambda(s)\to\lambda^{\mathcal{D}}$ as $\ s\to+\infty$, provided that  {\rm\textbf{(H1)}} and {\rm\textbf{(H2)}} hold.
 \end{proposition}

\begin{lemma}[Estimate of the upper limit]\label{upd}
Under the conditions of Proposition \ref{DT}, we have
 $$\limsup_{s\to+\infty}\lambda(s)\le \lambda^{\mathcal{D}}.$$
\end{lemma}
\begin{proof}
We construct a test function $\varphi_0(r)$ as follows: $\varphi_0(r)=0$ for $r\in[0,a)\cup(b,1]$ and $\varphi_0(r)=\varphi^{\mathcal{D}}(r)$ for $r\in[a,b]$ where $\varphi^{\mathcal{D}}(r)$ is the principal eigenfunction defined in Definition \ref{d2}.
 By the variational characterization of $\lambda(s)$ and $\lambda^{\mathcal{D}}$, and the fact $m(r)=0$ on $[a,b]$, we obtain
 \begin{eqnarray*}
 \lambda(s)\le\frac{\int_0^1 r^{d-1} e^{2sm(r)}\left(|\varphi_0'|^2+c(r)|\varphi_0|^2\right)\mathrm{d} r}{\int_0^1 r^{d-1} e^{2sm(r)}|\varphi_0|^2\mathrm{d} r}=
 \frac{\int_a^b r^{d-1} \left(|\varphi_0'|^2+c(r)|\varphi_0|^2\right)\mathrm{d} r}{\int_a^b r^{d-1} |\varphi_0|^2\mathrm{d} r}=
 \lambda^{\mathcal{D}}.
 \end{eqnarray*}
Taking the upper limit as $s\to+\infty$ on both sides of the above inequality, we complete the proof.
\end{proof}

Next, we address the boundary condition of the equation in Lemma \ref{converge}
\begin{equation*}\label{r}
-\varphi_*''(r)-\frac{d-1}{r}\varphi_*'(r)+c(r)\varphi_*(r)=\lambda_*\varphi_*(r)\quad r\in(a,b),
\end{equation*}
 which plays a crucial role in driving the lower bound of $\lambda(s)$. We will first establish the boundary condition of $\varphi_*(r)$ in the following lemma, which states that $\varphi_*(a)=\varphi_*(b)=0 $.\par

\begin{lemma}\label{id}
If $m(r)\in S_{\mathcal{D}}$, then the function $\varphi_*(r)$ derived from Lemma \ref{converge} subject to Dirichlet boundary condition, that is,
$$
 \varphi_*(a)=\varphi_*(b)=0.
$$
\end{lemma}
\begin{proof}
 Assume by contradiction that $\varphi_*(a)>0$. Let $s_i$ and $\varphi_1(s_i;r)$ be as in Lemma  \ref{converge}, and define
$$
N_i=\inf\left\{j\in \mathbb{N}\bigg|\varphi_1(s_{i};r)>\frac{1}{2}\varphi_*(a),\ r\in\bigcup_{n\ge j}[y_n,x_{n+1}]\right\}.
$$
Since
$$
\lim_{j\to+\infty}\inf_{r\in\bigcup_{n\ge j}[y_n,x_{n+1}]}\varphi_1(s_i;r)=\lim_{r\to a^-}\varphi_1(s_i;r)>\frac{1}{2}\varphi_*(a),
$$
for large $i$ by the convergence of $\varphi_1(s_i;r)$ as $i\to+\infty$ in Lemma \ref{converge}.
Thus, $N_i$ is well-defined.

In light of the definition of $\varphi_1(s_i;x)$ and \eqref{ubd}, we obtain
\begin{equation}\label{ky}
c^*\ge\lambda(s_i)\ge \int_\delta^a r^{d-1}e^{2s_im(r)}(|\varphi_1'|^2+c(r)|\varphi_1|^2)\mathrm{d} r,\quad\forall\ i>0,
\end{equation}
which, in particular, implies
\begin{equation}\label{key}
c^*\ge  \int_{y_{N_i}}^{x_{N_i+1}}r^{d-1}e^{2s_{i}m(r)} c(r)|\varphi_1|^2\mathrm{d} r\ge\frac{1}{4}\delta^{d-1}e^{2s_i\nu h^{N_i}}c_*\varphi_*^2(a)\beta^{N_i+1},
\end{equation}
due to $m(r)\in S_{\mathcal{D}}$.
Thanks to the fact that $s_i\to+\infty $ as $i\to+\infty$, along with \eqref{key}, we derive that
\begin{equation}\label{NN}
\lim_{i\to+\infty} N_i=+\infty.
\end{equation}

By the definition of $N_i$, there exists $Z\in[y_{N_i-1},x_{N_i})$ such that $\varphi(s_i;Z)\le\varphi_*(a)\mbox{/}2$. Therefore, we have
\begin{eqnarray}\label{kk}
\frac{1}{2}\varphi_*(a)
&\le&\varlimsup_{i\to+\infty}|\varphi_1(s_i;a)-\varphi_1(s_i;Z)|\nonumber\\
&\le&\varlimsup_{i\to+\infty}\left(|\varphi_1(s_i;Z)-\varphi_1(s_i;x_{N_i})|+\sum_{n=N_i}^{+\infty}|\varphi_1(s_i;x_n)-\varphi_1(s_i;y_n)|\nonumber\right.\\
&&\left.+\sum_{n=N_i}^{+\infty}|\varphi_1(s_i;y_n)-\varphi_1(s_i;x_{n+1})|\right).
\end{eqnarray}
We will drive a contradiction by showing that all three terms on the right hand side of \eqref{kk} are zero.

Applying Newton-Leibniz formula, we find
\begin{eqnarray*}\label{e11}
|\varphi_1(s_i;Z)-\varphi_1(s_i;x_{N_i})|&\le&\int_Z^{x_{N_i}}|\varphi_1'|\mathrm{d} r
\le\int_{y_{N_i-1}}^{x_{N_i}}|\varphi_1'|\mathrm{d} r\le \left(\int_{y_{N_i-1}}^{x_{N_i}}|\varphi_1'|^2\mathrm{d} r\right)^{\frac{1}{2}}\beta^{\frac{N_i}{2}}\nonumber\\
&\le&\left(\delta^{1-d}c^*\right)^{\frac{1}{2}}\beta^{\frac{N_i}{2}}
\end{eqnarray*}
where equation \eqref{ky} is used  in the last inequality and thanks to  \eqref{NN}, we have
\begin{eqnarray}\label{e11}
|\varphi_1(s_i;Z)-\varphi_1(s_i;x_{N_i})|\to 0,\quad i\to+\infty.
\end{eqnarray}
Similarly, for $n\in\mathbb{N}$, we can derive
\begin{eqnarray}\label{e21}
|\varphi_1(s_i;x_{n+1})-\varphi_1(s_i;y_n)|\le\left(\delta^{1-d}c^*\right)^{\frac{1}{2}}\beta^{\frac{n+1}{2}}.
\end{eqnarray}
Thus, again using \eqref{NN}, we get
\begin{eqnarray}\label{e22}
\sum_{n=N_i}^{+\infty}|\varphi_1(s_i;x_{n+1})-\varphi_1(s_i;y_n)|\le\left(\delta^{1-d}c^*\right)^{\frac{1}{2}}\sum_{n=N_i}^{+\infty}\beta^{\frac{n+1}{2}}\to 0,\quad \text{as } i\to+\infty.
\end{eqnarray}

Next, let us estimate the remainder term $\sum_{n=N_i}^{+\infty}|\varphi_1(s_i;x_n)-\varphi_1(s_i;y_n)|$. First, by applying Newton-Leibniz formula and H$\rm{\ddot{o}}$lder inequality, we obtain
\begin{eqnarray}\label{e31}
|\varphi_1(s_i;x_{n})-\varphi_1(s_i;y_n)|^2\le\alpha^{n+1}\int_{x_n}^{y_n}|\varphi_1'|^2\mathrm{d} r.
\end{eqnarray}
Utilizing \eqref{ky}, $m(r)\in S_{\mathcal{D}}$, and the above estimate \eqref{e31}, we get
\begin{eqnarray*}
c^*&\ge& \delta^{d-1}\sum_{n=N_{i}}^{+\infty}\left(\int_{x_n}^{y_n}e^{2s_{i}\tilde{m}}|\varphi_1'|^2\mathrm{d} r+c_*\int_{y_n}^{x_{n+1}}e^{2s_{i}\tilde{m}}|\varphi_1|^2\mathrm{d} r\right)\nonumber\\
&=&\delta^{d-1}\sum_{n=N_{i}}^{+\infty}\left(e^{-2s_{i}h^{n}}\int_{x_n}^{y_n}|\varphi_1'|^2\mathrm{d} r
+c_*e^{2s_{i}\nu h^{n}}\int_{y_n}^{x_{n+1}}|\varphi_1|^2\mathrm{d} r\right)\nonumber\\
&\ge&\delta^{d-1}\sum_{n=N_{i}}^{+\infty}\left(e^{-2s_{i}h^{n}}\frac{|\varphi_1(s_{i};y_n)-\varphi_1(s_{i};x_n)|^2}{\alpha^{n+1}}
+\frac{1}{4}\varphi_*^2(a) c_*e^{2s_{i}\nu h^{n}}\beta^{n+1}\right)\nonumber\\
&\ge&\delta^{d-1}\sqrt{c_*}\varphi_*(a)\inf_{n\ge N_{i}}\left(\frac{\beta}{\alpha}\right)^{\frac{n+1}{2}}\sum_{n=N_{i}}^{+\infty}|\varphi_1(s_{i};y_n)-\varphi_1(s_{i};x_n)|.
\end{eqnarray*}
Thus, by the fact $\beta>\alpha$, we derive that
\begin{equation}\label{DDE}
\sum_{n=N_i}^{+\infty}|\varphi_1(s_i;x_n)-\varphi_1(s_i;y_n)|\le \frac{c^*}{\delta^{d-1}\sqrt{c_*}\varphi_*(a)}\left(\frac{\alpha}{\beta}\right)^{\frac{N_i+1}{2}}\to 0, \quad \text{as } i\to+\infty.
\end{equation}

Thus, combining \eqref{kk}, \eqref{e11}, \eqref{e22}, and \eqref{DDE} yields that $\frac{1}{2}\varphi_*(a)\le 0$, which is absurd. Hence, we deduce that $\varphi_*(a)=0$ and, consequently, $\varphi_*(b)=0$ by the similar manner.
 \end{proof}

In the following, we give the lower bound of $\liminf_{s\to+\infty}\lambda(s)$.
\begin{lemma}[Estimate of the lower limit]\label{lowerd}
Under the conditions of Proposition \ref{DT}, we have
 $$
\liminf_{s\to+\infty}\lambda(s)\ge\lambda^{\mathcal{D}}.
 $$
 \end{lemma}

\begin{proof}
We consider the variational characterization of the principal eigenvalue in \eqref{voa} with the transformation $w_1=e^{sm(r)}r^{\frac{d-1}{2}}\varphi_1$, where the sequence $\{s_i \}_{i=1}^{+\infty}$ is the one obtained in the Lemma \ref{converge}. Recall that $m(r)=0$ on $[a,b]$, which implies $w_1=r^{d-1}\varphi_1$ on $[a,b]$. Then we obtain 
\begin{align}
&\lambda(s_{i})\\
=&\int_0^1\left(\left|w_1'(s_i;r)-sm'(r)w_1(s_i;r)-\frac{d-1}{2r}w_1(s_i;r)\right|^2+c(r)|w_1(s_{i};r)|^2\right)\mathrm{d} r\nonumber\\
\ge&\int_0^ac(r)|w_1(s_{i};r)|^2\mathrm{d} r
+\!\!\int_a^b r^{d-1}(|\varphi_1'(s_{i};r)|^2+c(r)|\varphi_1(s_{i};r)|^2)\mathrm{d} r
+\!\!\int_b^1c(r)|w_1(s_{i};r)|^2\mathrm{d} r
\nonumber\\
\triangleq& I_1+I_2+I_3.
\label{estd'}
\end{align}
For any $\epsilon>0$, we will establish the following three estimates.
\begin{eqnarray}
\liminf_{i\to+\infty}I_1&\ge&\lambda^{\mathcal{D}}\mu([0,a+3\epsilon])-C\epsilon\label{1},\\
\liminf_{i\to+\infty}I_2&\ge&\lambda^{\mathcal{D}}\mu([a+2\epsilon,b-2\epsilon])\label{2},\\
\liminf_{i\to+\infty}I_3&\ge&\lambda^{\mathcal{D}}\mu([b-3\epsilon,1])-C\epsilon,\label{3}
\end{eqnarray}
where $C=C(d,c^*,a,b)$ is some positive constant.
Once the above estimates are established, we can deduce that
\begin{equation}\label{x}
\liminf_{i\to+\infty}\lambda(s_i)\ge\lambda^{\mathcal{D}}\mu([0,1])-2C\epsilon.
\end{equation}
Therefore, we derive that $\liminf_{i\to+\infty}\lambda(s_i)\ge\lambda^{\mathcal{D}}$ by sending $\epsilon\to 0^+$ in \eqref{x}.

Recall that $\int_0^1w^2(s_i;r)\mathrm{d} r=1$ for each $i>0$, which implies $\{w^2(s_i,\cdot)\}_{i>0}$ is weakly compact in $L^1((0,1))$. Thus, we could select a subsequence $\{s_i\}_{i=1}^{+\infty}$, denoted by itself, satisfying
\begin{equation}\label{weakc}
\lim_{i\to +\infty}\int_0^1w^2(s_i;r)\zeta(r)\mathrm{d} r=\int_{[0,1]}\zeta(r)\mu(\mathrm{d} r), \quad \forall\ \zeta(r)\in C([0,1]),
\end{equation}
where $\mu$ is a certain probability measure.

To estimate \eqref{1}, we define
 \begin{equation*}\label{z1}
\xi(r) = \left\{ \begin{array}{lll}
c(r), &r\in [0,a)\\
c(a), &r\in (a, a+3\epsilon]\\
-\frac{c(a)}{\epsilon}(r-a-4\epsilon),& (a+3\epsilon,a+4\epsilon]\\
0, & (a+4\epsilon,1].
 \end{array} \right.
\end{equation*}
From \eqref{weakc}, we obtain
\begin{eqnarray}
\liminf_{i\to+\infty}I_1
&\ge&\liminf_{i\to+\infty}\int_{0}^{1}w_1^2(s_{i};r)\xi(r)\mathrm{d} r
-\limsup_{i\to+\infty}\int_{a}^{a+4\epsilon}w_1^2(s_{i};r)\xi(r)\mathrm{d} r\nonumber\\
&=&\int_{[0,1]}\xi(r)\mu(\mathrm{d} r)-\limsup_{i\to+\infty}\int_{a}^{a+4\epsilon}r^{d-1}\varphi_1^2(s_{i};r)\xi(r)\mathrm{d} r\nonumber\\
&\ge& \int_{[0,a+3\epsilon]}\xi(r)\mu(\mathrm{d} r)-C\epsilon
\ge \lambda^{\mathcal{D}}\mu([0,a+3\epsilon])-C\epsilon,\label{de'}
\end{eqnarray}
where the second inequality follows from the boundedness of $\varphi_1(s_i;r)$ as given in \eqref{wb}, and the last inequality follows from \textbf{(H2)}: $c(r)>\lambda^{\mathcal{D}}$ on $[0,a]\cup[b,1]$.  And
in a similar manner, we can obtain the estimate \eqref{3}.

To estimate \eqref{2}, we will consider it in two cases.

\textbf{Case 1}. $
\liminf\limits_{i\to+\infty}\int_a^br^{d-1}\varphi_1^2(s_{i};r)\mathrm{d} r=0.
$

Set the  continuous function $\xi(r)=1$ on $[a+\epsilon,b-\epsilon]$, $\xi(r)=0$ on $[o,a]\cup[b,1]$, and $0\le \xi \le 1$ on $[0,1]$.
 Then, according to equation \eqref{weakc}, we obtain
\begin{align}\label{same}
\liminf_{{i\to+\infty}}\int_{a}^{b}r^{d-1}\varphi_1^2(s_{i};r)\mathrm{d} r
&=\liminf_{{i\to+\infty}}\int_{a}^{b}w_1^2(s_{i};r)\mathrm{d} r\nonumber\\
&\ge\liminf_{{i\to+\infty}}\int_{0}^{1}w_1^2(s_{i};r)\xi(r)\mathrm{d} r\nonumber\\
&=\int_{[0,1]}\xi(r)\mu(\mathrm{d} r)\nonumber\\
&\ge\int_{[a+\epsilon,b-\epsilon]}\xi(r)\mu(\mathrm{d} r)\nonumber\\&\ge\mu([a+2\epsilon,b-2\epsilon]).
\end{align}
Thus, we have that $\mu([a+2\epsilon,b-2\epsilon])=0$ and $I_2\ge 0=\lambda^{\mathcal{D}}\mu([a+2\epsilon,b-2\epsilon])$.

\textbf{Case 2}. $
\liminf\limits_{i\to+\infty}\int_a^br^{d-1}\varphi_1^2(s_{i};r)\mathrm{d} r>0.
$

We  denote
$$
v(s_{i};r)\triangleq\frac{\varphi_1(s_{i};r)}{\sqrt{\int_a^br^{d-1}|\varphi_1(s_{i};r)|^2\mathrm{d} r}}.
$$
According to Lemma \ref{converge}, we can select a subsequence of $\{s_i\}$, still denoted by $\{s_i\}$, such that
\begin{equation}\label{w}
v(s_{i};r)\rightharpoonup v_1(r), \quad \text{in }H^1((a,b)) \text{  as }i\to+\infty,
\end{equation}
and
\begin{equation}\label{5}
v(s_{i};r)\to \frac{\varphi_*(r)}{\sqrt{\int_a^br^{d-1}|\varphi_*(r)|^2\mathrm{d} r}}\triangleq v_1(r),\ \text{in } C([a,b])\cap C^1_{\mathrm{loc}}((a,b)) \text{  as }i\to+\infty,
\end{equation}
with $\int_a^br^{d-1}v_1^2(r)\mathrm{d} r=1$. From Lemma \ref{id}, we have $\varphi_*(a)=\varphi_*(b)=0$, and it follows that
 $v_1(a)=v_1(b)=0$ .\par
Now we derive
\begin{eqnarray*}
\liminf_{i\to+\infty}I_2
&=&\liminf_{i\to+\infty}\left(\int_a^br^{d-1}|\varphi_1(s_{i};r)|^2\mathrm{d} r\int_a^br^{d-1}\left(|v'(s_{i};r)|^2+c(r)|v(s_{i};r)|^2\right)\mathrm{d} r\right)\nonumber\\
&\ge&\mu([a+2\epsilon,b-2\epsilon])\liminf_{i\to+\infty}\int_a^br^{d-1}\left(|v'(s_{i};r)|^2+c(r)|v(s_{i};r)|^2\right)\mathrm{d} r\nonumber\\
&\ge&\mu([a+2\epsilon,b-2\epsilon])\int_a^br^{d-1}\left(|v_1'(r)|^2+c(r)|v_1(r)|^2\right)\mathrm{d} r\nonumber\\
&\ge&\mu([a+2\epsilon,b-2\epsilon])\lambda^{\mathcal{D}},\label{ddab}
\end{eqnarray*}
where we have used \eqref{same} in the first inequality, the weakly lower semi-continuity of the $L^2$ norm \eqref{w} and \eqref{5} in the second inequality, and the variational characterization of $\lambda^{\mathcal{D}}$ in the last inequality.
\end{proof}

\subsection{Consider $m(r)\in S_{\mathcal{N}}$}\label{DD}
Similarly to the case where $m(r)\in S_{\mathcal{D}}$, we obtain the following convergence result for $m(r)\in S_{\mathcal{N}}$.
\begin{proposition}\label{NT}
 If $m(r)\in S_{\mathcal{N}}$,  then the principal eigenvalue $\lambda(s)$ of \eqref{eq} satisfies
 $$\lim_{s\to+\infty}\lambda(s)=\lambda^{\mathcal{N}},$$ provided that  {\rm\textbf{(H1)}} and {\rm\textbf{(H2)}} hold.
 \end{proposition}

\begin{lemma}[Estimate for the lower limit]\label{nnl}
Under the conditions of Proposition \ref{NT}, we have
$$\liminf_{s\to+\infty}\lambda(s)\ge \lambda^{\mathcal{N}}.$$
\end{lemma}

\begin{proof}
We recall that  $\varphi_1(s;r)$ is the principal eigenfunction of \eqref{eq} with respect to $\lambda(s)$ for any fixed $s>0$.
By the variational characterization \eqref{voe} and the transformation $w_1(s_{i};r)=r^{\frac{d-1}{2}}e^{s_im(r)}\varphi_1(s_{i};r)$,
we obtain
\begin{align}\label{ne}
\lambda(s_{i})=&\int_0^1r^{d-1}e^{2s_im(r)}(|\varphi_1'(s_{i};r)|^2+c(r)|\varphi_1(s_{i};r)|^2)\mathrm{d} r\nonumber\\
\ge&\left(\int_{0}^a+\int_b^{1}\right)c(r)w_1^2(s_{i};r)dr
+\int_a^br^{d-1} (|\varphi_1'(s_{i};r)|^2+c(r)|\varphi_1(s_{i};r)|^2)\mathrm{d} r,
\end{align}
where the sequence $\{s_i\}_{i=1}^{\infty}$ is the one obtained in Lemma \ref{converge}.
In light of \eqref{DDN}, using the estimate in \eqref{1} and \eqref{3}, we have
\begin{align}\label{nne}
\liminf_{i\to+\infty}\left(\int_{0}^a+\int_b^{1}\right)c(r)w_1^2(s_{i};r)\mathrm{d} r\ge&\lambda^{\mathcal{D}}\left[\mu([0,a+3\epsilon])+\mu([b-3\epsilon,1])\right]-C\epsilon\nonumber\\
>&\lambda^{\mathcal{N}}\left[\mu([0,a+3\epsilon])+\mu([b-3\epsilon,1])\right]-C\epsilon.
\end{align}
For the last term in \eqref{ne}, we claim that
\begin{equation}\label{abnn}
\liminf_{i\to+\infty}\int_a^b r^{d-1} (|\varphi_1'(s_{i};r)|^2+c(r)|\varphi_1(s_{i};r)|^2)\mathrm{d} r
\ge\lambda^{\mathcal{N}}\mu([a+2\epsilon,b-2\epsilon]).
\end{equation}
Assuming \eqref{abnn} holds and combining the estimates from \eqref{nne} and \eqref{abnn}, we obtain the desired result $\liminf_{s\to+\infty}\lambda(s)\ge \lambda^{\mathcal{N}}$ by first taking $i\to\infty$ and then letting $\epsilon \to 0^+$ on the both sides of \eqref{ne}.

Next, we prove our claim \eqref{abnn} in two cases.
\\
\textbf{Case 1}. $
\liminf\limits_{i\to+\infty}\int_a^br^{d-1}\varphi_1^2(s_{i};r)\mathrm{d} r=0.
$

By \eqref{same} in Section 3, we have
 $\mu([a+2\epsilon,b-2\epsilon])=0$, which implies that \eqref{abnn} holds.
\\
\textbf{Case 2}.
$
\liminf\limits_{i\to+\infty}\int_a^br^{d-1}\varphi_1^2(s_{i};r)\mathrm{d} r>0.
$

Set
$$
v(s_{i};r)=\frac{\varphi_1(s_{i};r)}{\sqrt{\int_a^br^{d-1}\varphi_1^2(s_{i};r)\mathrm{d} r}}\in H^1((a,b)).
$$
Relying on the property that $\int_a^br^{d-1}v^2(s_i;r)\mathrm{d} r=1$, \eqref{same} and the variational characterization of $\lambda^{\mathcal{N}}$, we obtain
\begin{eqnarray*}
&&\liminf_{i\to+\infty}\int_a^br^{d-1}(|\varphi_1'(s_{i};r)|^2+c(r)|\varphi_1(s_{i};r)|^2)\mathrm{d} r\nonumber\\
&=&\liminf_{i\to+\infty}\int_a^br^{d-1}\varphi_1^2(s_{i};r)\mathrm{d} r\int_a^br^{d-1} (|v'(s_{i};r)|^2+c(r)|v(s_{i};r)|^2)\mathrm{d} r\nonumber\\
&\ge&\lambda^{\mathcal{N}}\mu([a+2\epsilon,b-2\epsilon]).
\end{eqnarray*}
The proof is completed.
\end{proof}

Next, we will establish the upper bound. Before doing so, we introduce the following definitions.
\begin{definition}\label{dnn}
Let $m(r)\in S_{\mathcal{N}}$ and let, $\delta$, $\alpha$, $\beta$, $h$, $\nu$, $l=0$ be as defined in Definition \ref{D3}.
We define
\begin{equation*}
\sigma_n(s)\triangleq(\alpha\beta)^{\frac{n}{2}} e^{s (1+\nu) h^{n} },
\end{equation*}
and
\begin{equation*}
 p_n(s)\triangleq
\frac { 1 } { \prod _ { k = n } ^ {+ \infty } \left( 1 + (\alpha\beta)^{\frac{k}{2}} e^{ s (1+\nu)h^{k}} \right) }=\frac { 1 } { \prod _ { k = n } ^ {+ \infty } \left( 1 + \sigma_k(s) \right) }.
\end{equation*}
\end{definition}

\begin{remark}\label{propnn}
From Definition \ref{dnn}, it is straightforward to check the following properties on $p_n(s)$. \par
{\rm(1)} $p_n(s)$ is increasing in $n$, $\lim_{n\to +\infty}p_n(s)=1$ for any fixed $s>0$ and $\lim_{s\to +\infty}p_n(s)=0$ for any fixed $n\in \mathbb{N}$;\par
{\rm(2)} $p_{n+1}(s)-p_n(s)=\sigma_n(s)p_n(s)$.
\end{remark}

Now, we can proceed to deduce the upper bound estimate.

\begin{lemma}[Estimate for the upper limit]\label{nnu}
Under the conditions of Proposition \ref{NT}, we have
$$\limsup_{s\to+\infty}\lambda(s)\le \lambda^{\mathcal{N}}.$$
\end{lemma}

\begin{proof}
The condition $m(r)\in S_{\mathcal{N}}$  ensures the existence of constants $\delta$, $\alpha$, $\beta$, $h$, $\nu$ such that $m(r)\leq\bar{m}(r)$ in $[\delta,a)\cup(b,1-\delta]$. Recall that $\bar{m}(\delta)=h$, so by the continuity of $m(r)$, there exists a small positive constant $\delta_1$ such that $m(r)\leq \frac{3}{2}h$ for $r\in (\delta-\delta_1,\delta]$.
We define a test function $\varphi_0\in H^1([0,1])$ as follows.
\begin{equation}\label{var1}
\varphi _ { 0 }(s;r) = \left\{ \begin{array}{lcl}
0, && r\in[0,\delta-\delta_1),\\
\frac{\varphi^{\mathcal{N}}(a)p_1(s)}{\delta_1}(r-\delta+\delta_1), && r\in[\delta-\delta_1,\delta),\\
\varphi^{\mathcal{N}}(a)p_n(s), & & r\in [Y_{n-1}, X_{n}),\\
\varphi^{\mathcal{N}}(a)\left(p_n(s)+\frac { p_{n+1}(s) - p _ { n }(s) } { \beta ^{ n } }(r-X_n)\right), & & r\in [X_{n}, Y_n), \,\, \\
\varphi^{\mathcal{N}}(r),&& r\in[a,b], \end{array} \right.
\end{equation}
and for $r\in(b,1]$, $
\varphi_0(s;r)=\frac{\varphi^{\mathcal{N}}(b)}{\varphi^{\mathcal{N}}(a)} \varphi_0(s;1-r)$,
where the function
$\varphi^{\mathcal{N}}(r)$ is the principal eigenfunction corresponding to $\lambda^{\mathcal{N}}$ defined in Definition \ref{d2}, the sequences $\{X_n\}_{n=1}^{\infty}$, $\{Y_{n-1}\}_{n=1}^{\infty}$ are defined in Definition \ref{D3}, and the function $p_n(s)$ is given in Definition \ref{dnn}. Then, by the variational characterization of $\lambda(s)$ in \eqref{voa}, we have
\begin{eqnarray}
\lambda(s)
&\le&\frac{ \left(\int_0^a+\int_a^b+\int_b^1\right)r^{d-1}e^{2sm(r)}\left(|\varphi_0'(s;r)|^2+c(r)|\varphi_0(s;r)|^2\right)\mathrm{d} r} {\int_0^1r^{d-1}|\varphi_0(s;r)|^2\mathrm{d} r}\nonumber\\
&\le& \lambda^{\mathcal{N}}+ \frac{ I_0^a[\varphi_0](s)+I_b^1[\varphi_0](s)}{\int_a^br^{d-1}|\varphi^{\mathcal{N}}(r)|^2\mathrm{d} r},\label{lbd1}
\end{eqnarray}
where we denote the notation $$I_{e_1}^{e_2}[\varphi](s)\triangleq\int_{e_1}^{e_2}r^{d-1}e^{2sm(r)}\left(|\varphi'(r)|^2+c(r)|\varphi(r)|^2\right)\mathrm{d} r,$$
for any function $\varphi(r)\in H^1([0,1])$ and $e_1<e_2$ in $(0,1)$.

\par
We begin by considering $I_0^a[\varphi_0](s)$. Due to the symmetry of $\varphi_0(s;r)$, $I_b^1[\varphi_0](s)$ can be treated similarly. We have
\begin{align}
&I_0^a[\varphi_0](s)\nonumber\\
=& \int_{0}^{\delta}r^{d-1}e^{2sm(r)}\left(|\varphi_0'|^2+c(r)|\varphi_0|^2\right)\mathrm{d} r+\int_\delta^a r^{d-1}e^{2sm(r)}\left(|\varphi_0'|^2+c(r)|\varphi_0|^2\right)\mathrm{d} r\nonumber\\
\le&\int_{\delta-\delta_1}^{\delta} r^{d-1} e^{2sm(r)}\left(|\varphi_0'|^2+c(r)|\varphi_0|^2\right)\mathrm{d} r+\sum_{n=1}^{+\infty} \int_{Y_{n-1}}^{Y_{n}}r^{d-1}e^{2s\bar{m}(x)}\left(|\varphi_0'|^2+c(r)|\varphi_0|^2\right)dx\nonumber\\
\le&C\left[e^{3sh}p_1^2(s)+ \sum_{n=1}^{+\infty} \alpha^{n} e^{2sh^{n}}p_n(s)^2
+\sum_{n=1}^{+\infty}\left(e^{-2s\nu h^{n}} \frac{(p_{n+1}(s)-p_n(s))^2}{\beta^{n}}+\beta^{n}e^{-2s\nu h^{n}}p_{n+1}^2(s)\right)\right]\nonumber\\
\le&C\left(e^{3sh}p_1^2(s)+\sum_{n=1}^{+\infty}\alpha^{n}e^{2sh^{n}}p_n^2(s)+\sum_{n=1}^{+\infty}\beta^{n}e^{-2s\nu h^{n}}p_{n+1}^2(s)\right)\nonumber\\
\triangleq& C(E(s)+F(s)+G(s)),\label{est E(s) 0}
\end{align}
where $C=C(\varphi^{\mathcal{N}}(a), a, d,  c^*)$ is a positive constant, and the last inequality follows from (2) in Remark \ref{propnn}.
To study the upper bound of \eqref{est E(s) 0}, we estimate these three terms on the righthand side of it separately.
\par
\textbf{Estimate of $E(s)$.}
Under the fact that $\nu>1$, we can derive
\begin{equation}
E(s)=\frac { e^{3sh} } { \prod _ { k = 1 } ^ {+ \infty } \left( 1 + (\alpha\beta)^{\frac{k}{2}} e^{s (1+\nu) h^{k} } \right)^2 }\le\frac { e^{3sh} } { \alpha\beta e^{2s (1+\nu) h }  } \to 0, \quad \text{ as } s\to+\infty.
\end{equation}

\textbf{Estimate of $G(s)$.} Since $\beta^{n}e^{-2s\nu h^{n}}p_{n+1}^2\le \beta^{n}$,  we have
\begin{equation}\label{les}
\lim_{s\to+\infty}G(s)=\lim_{s\to+\infty}\sum_{n=1}^{+\infty}\beta^{n}e^{-2s\nu h^{n}}p_{n+1}^2(s)=\sum_{n=1}^{+\infty}\lim_{s\to+\infty}\beta^{n}e^{-2s\nu h^{n}}p_{n+1}^2(s)=0,
\end{equation}
due to the uniformly convergence theorem.

\textbf{Estimate of $F(s)$.}
For any $\epsilon\in (0, 1)$, using the definition of $\sigma_n(s)$ in Definition \ref{dnn}, we can find
$0<K_1(s,\epsilon)\le K_2(s,\epsilon)$ such that
\begin{equation}\label{K1K2}
\left\{ \begin{array} { l } \sigma_1(s),\ \sigma_2(s),...\ \sigma_{K_1}(s)\in (\frac{1}{\epsilon},+\infty), \\
\sigma_{K_1+1}(s),\ \sigma_{K_1+2}(s),...\ \sigma_{K_2}(s)\in (\epsilon, \frac{1}{\epsilon}),\\
\sigma_{K_2+1}(s),\ \sigma_{K_2+2}(s),...\in (0,\epsilon), \end{array} \right.
\end{equation}
where we use the notations $K_1$ and $K_2$ for simplicity. Note that for any fixed $\epsilon$, $K_1$ and $K_2$ are functions depending on $s$. Based on this, we insert $p_n(s)$ into $F(s)$ and divide it into three parts.
\begin{eqnarray}\label{DF}
F(s)
=\left(\sum_{n=1}^{K_1}+\sum_{n=K_1+1}^{K_2}+\sum_{n=K_2+1}^{+\infty}\right)\frac {\alpha^{n}e^{2sh^{n}}} { \prod _ { k = n } ^ { +\infty } \left( 1 + \sigma _ { k }(s) \right)^2 }
\triangleq F_1(s)+F_2(s)+F_3(s).
\end{eqnarray}
Next, we estimate $F_1(s)$, $F_2(s)$ and $F_3(s)$.\par

For $F_1(s)$, recall that $0<\alpha<\beta<1$, $\nu>1$. We compute
\begin{equation}\label{u1}
\sigma_n(s)=(\alpha\beta)^{\frac{n}{2}} e^{s (1+\nu) h^{n} }\ge \alpha^{n}e^{2sh^{n}},
\end{equation}
and
\begin{eqnarray}\label{u2}
\sum_{n=1}^{K_1}\frac{\sigma_n}{ \prod _ { k = n } ^ { +\infty } \left( 1 + \sigma _ { k }(s) \right)}&=&\sum_{n=1}^{K_1}\left(\frac{1}{ \prod _ { k = n+1 } ^ { +\infty } \left( 1 + \sigma _ { k }(s) \right)}-\frac{1}{ \prod _ { k = n } ^ { +\infty } \left( 1 + \sigma _ { k }(s) \right)}\right)\nonumber\\
&\le&\frac{1}{ \prod _ { k = K_1+1 } ^ { +\infty } \left( 1 + \sigma _ { k }(s) \right)}\le 1.
\end{eqnarray}
In view of the definition of $K_1$, combining \eqref{u1} and \eqref{u2} yields that
\begin{eqnarray}
F_1(s)
&=&\sum_{n=1}^{K_1}\frac{\alpha^{n}e^{2sh^{n}}}{ \prod _ { k = n } ^ { +\infty } \left( 1 + \sigma _ { k }(s) \right)^2}
\le\sum_{n=1}^{K_1}\frac{\sigma_n}{ \prod _ { k = n } ^ { +\infty } \left( 1 + \sigma _ { k }(s) \right)}\cdot \frac{1}{ \prod _ { k = n } ^ { +\infty } \left( 1 + \sigma _ { k }(s) \right)}\nonumber\\
&\le&\sum_{n=1}^{K_1}\frac{\sigma_n}{ \prod _ { k = n } ^ { +\infty } \left( 1 + \sigma _ { k }(s) \right)}\cdot \frac{1}{  \sigma _ { K_1 }(s) }
\le \epsilon\label{F1}.
\end{eqnarray}
 \par

In order to estimate $F_2(s)$, we
first show that the interval $(\epsilon, \frac{1}{\epsilon})$ contains at most one term for sufficiently large $s$. Assume by contradiction that there exist $s_j\to+\infty$ as $j\to+\infty$ such that $\sigma_K,\ \sigma_{K+1}\in (\epsilon, \frac{1}{\epsilon})$ with $K=K_1(s_j,\epsilon)+1$. Then we have
\begin{equation}
\epsilon<(\alpha\beta)^{\frac{K+1}{2}}e^{s_jh^{K+1}(1+\nu)}<(\alpha\beta)^{\frac{K}{2}}e^{s_jh^{K}(1+\nu)}<\frac{1}{\epsilon} \label{one}
\end{equation}
for large $j$.
The last inequality of \eqref{one} implies that
\begin{equation}
s_j<\frac{1}{h^{K}(1+\nu) }\ln\frac{1}{\epsilon(\alpha\beta)^{\frac{K}{2}}}.\label{s1}
\end{equation}
This means that $K\to+\infty$ as $j\to+\infty$. Substituting \eqref{s1} into the first inequality of \eqref{one}, we obtain
\begin{eqnarray}
\epsilon<(\alpha\beta)^{\frac{K+1}{2}}e^{s_jh^{K+1}(1+\nu)}<\frac{(\sqrt{\alpha\beta})^{K(1-h)+1}}{\epsilon^h}.\label{con1}
\end{eqnarray}
Since $\alpha, \beta, h\in(0,1)$, the right hand side of \eqref{con1} tends to zero as $j\to+\infty$, which leads to a contradiction. Therefore, the only possible term $\sigma_K\in(\epsilon,\frac{1}{\epsilon})$ for large $s$. Then we have
\begin{eqnarray}
F_2(s)&\leq&\frac{\alpha^{K}e^{2sh^{K}}}{ \prod _ { k = K } ^ { +\infty } \left( 1 + \sigma _ { k }(s) \right)^2}
\le\left(\frac{\alpha}{\beta}\right)^{\frac{K}{2}}\sigma_Ke^{s(1-\nu)h^K}\le\frac{1}{\epsilon}\left(\frac{\alpha}{\beta}\right)^{\frac{K}{2}},\label{F2}
\end{eqnarray}
where the last inequality follows from $\sigma_K\in(\epsilon,\frac{1}{\epsilon})$ and the fact that $\nu<1$.
\par

Finally, for $F_3(s)$, when $n\ge K_2+1$, we have $\sigma_n\in(0,\epsilon)$ from \eqref{K1K2}, and thus
\begin{equation*}\label{k2}
\sigma_n=(\alpha\beta)^{\frac{n}{2}} e^{s (1+\nu) h^{n} }\le (\alpha \beta)^{\frac{n-K_2-1}{2}}\sigma_{K_2+1}\le(\alpha \beta)^{\frac{n-K_2-1}{2}}\epsilon.
\end{equation*}
Hence, by \eqref{u1}, we obtain
\begin{eqnarray}
F_3(s)
=\sum_{n=K_2+1}^{+\infty}\frac{\alpha^{n}e^{2sh^{n}}}{ \prod _ { k = n } ^ { +\infty } \left( 1 + \sigma _ { k }(s) \right)^2}
\le\sum_{n=K_2+1}^{+\infty}\sigma_n
\le\epsilon\sum_{n=K_2+1}^{+\infty}(\alpha \beta)^{\frac{n-K_2-1}{2}}
\le C\epsilon,\label{F3}
\end{eqnarray}
where the last inequality follows from the fact $0<\alpha$, $\beta<1$. Therefore, combining \eqref{F1}, \eqref{F2} with \eqref{F3}, we get the estimate for $F(s)$ for large $s$
\begin{eqnarray}\label{F}
F(s)=F_1(s)+F_2(s)+F_3(s)
\le\epsilon+\frac{1}{\epsilon}\left(\frac{\alpha}{\beta}\right)^{\frac{K}{2}}+C\epsilon.
\end{eqnarray}
Recall that $K\to\infty$ as $s\to\infty$. Sending $s\to +\infty$, and then letting $\epsilon\to 0^+$, we conclude from \eqref{F} that $F(s)\to 0$ as $s\to+\infty$.\par

\par
By the estimates of $E(s)$, $F(s)$ and $G(s)$, we conclude $I_0^a[\varphi_0](s)\to0$ as $s\to+\infty$. Similarly, we can show that $I_b^1[\varphi_0](s)\to0$ as $s\to+\infty$. Taking the upper limit on both sides of \eqref{lbd1} completes the proof.
\end{proof}

\textbf{Proof of Proposition \ref{NT}.} Combining Lemma \ref{nnl} and Lemma \ref{nnu}, we obtain Proposition \ref{NT}.

\section{Proof of Theorem \ref{robin}}\label{notconver}

In order to establish Theorem \ref{robin}, we begin by introducing a lemma that illustrates the continuity dependence of the principal eigenvalue $\lambda(s)$ on $m(r)$ for fixed $s$.  To emphasize the dependence of $\lambda(s)$ on the advection term $m(r)$, we write $\lambda(s)$ as $\lambda(s,m(r))$.

\begin{lemma}\label{lemma}
The principal eigenvalue
$\lambda(s,m(r))$ continuously depends on $m(r)\in  C([0,1])$.  Specifically, we have
\begin{equation}
\left|\lambda(s,m_1(r))- \lambda(s,m_2(r)) \right|\le c^*\left(e^{4s\|m_1-{m_2}\|_{C([0,1])}}-1\right),\label{rbes}
\end{equation}
where $c^*=\max_{[0,1]}c(r)$.
\end{lemma}

\begin{proof}
 By the variational characterization \eqref{voa}, we have
\begin{align}
&\lambda(s,m_1(r))- \lambda(s,m_2(r)) \nonumber\\
=&\min_{\int_0^1 r^{d-1} e^{2sm_1(r)}|\phi|^2\mathrm{d} r=1}\int_0^1 r^{d-1} e^{2sm_1(r)}\left[|\phi'|^2+c(r)|\phi|^2\right]\mathrm{d} r-\lambda(s,m_2(r)),\label{vart}
\end{align}
for all $\phi\in H^1([0,1])$.
Taking the test function $\phi(r)=\varphi_2(r)$, which is the principal eigenfunction corresponding to \eqref{eq} with  $m(r)=m_2(r)$, from \eqref{vart}, we obtain
\begin{align}
&\lambda(s,m_1(r))- \lambda(s,m_2(r))\nonumber\\
\le&\frac{\int_0^1 r^{d-1} e^{2sm_1(r)}\left[|\varphi_2'(r)|^2+c(r)|\varphi_2(r)|^2\right]\mathrm{d} r}{\int_0^1r^{d-1}e^{2sm_1(r)}|\varphi_2(r)|^2\mathrm{d} r}
-\lambda(s,m_2(r))\nonumber\\
\le& e^{4s\|m_1-{m_2}\|_{C([0,1])}}\frac{\int_0^1 r^{d-1} e^{2sm_2(r)}\left(|\varphi_2'(r)|^2+c(r)|\varphi_2(r)|^2\right)\mathrm{d} r}{\int_0^1 r^{d-1} e^{2sm_2(r)}|\varphi_2(r)|^2\mathrm{d} r}-\lambda(s,m_2(r))\nonumber\\
=& \lambda(s,m_2(r))\left(e^{4s\|m_1-{m_2}\|_{C([0,1])}}-1\right)
\le c^*\left(e^{4s\|m_1-{m_2}\|_{C([0,1])}}-1\right),\label{p}
\end{align}
where the last inequality follows from \eqref{ubd}. We estimate $\lambda(s,m_2(r))- \lambda(s,m_1(r))$ by using the same approach and thus obtain \eqref{rbes}.
\end{proof}
We are now in a position to prove Theorem \ref{robin} by using Proposition \ref{DT}, Proposition \ref{NT} and Lemma \ref{lemma}.
\par

\textbf{Proof of Theorem \ref{robin}.}
We aim to construct a symmetric counterexample, denoted as $m^*(r)$. In the following, we will construct it on the interval $[0,a]$, and the definition of $m^*(r)$ on the interval $[b,1]$ can be obtained immediately through symmetry.
To facilitate understanding of the construction process, we will divide it into three steps.

\textbf{Step 1.} Construction of a sequence of functions $\{m_n(r)\}_{n=1}^{+\infty}$.

Define $m_1(r)=m_{\mathcal{D}}(r)$, where $m_{\mathcal{D}}(r)$ is in Example \ref{EE}.
Recall that $m_{\mathcal{D}}(r)\in S_{\mathcal{D}}$, satisfying
$m_1(r)\ge 0$ on $[x_n,y_n]$ with $z_n=\frac{x_n+y_n}{2}$ and $m_1(z_n)=m_1'(z_n)=0$.
It is easy to see that
$$\lim_{s\to+\infty}\lambda(s,m_1(r))=\lambda^{\mathcal{D}},$$
from Theorem \ref{DT}.
 Thus, there exists $s_1>8\ln^{-1}(1+\frac{1}{2c^*})$ such that
 $$\big|\lambda(s_1,m_1(r))-\lambda^{\mathcal{D}}\big|<\frac{1}{2}.$$
 Let $\tau_1=s_1^{-2}$, thanks to $\lim_{r\to a^-}m_1(r)=0$, we can choose $\delta(\tau_1)>0$, such that $|m_1(r)|<\tau_1$, as $r\in(a-\delta(\tau_1),a)$. \par
We select a large $j_1\in \mathbb{N}$ such that $z_{j_1}>a-\delta(\tau_1)$, and define
$$
m_2(r)=\left\{
\begin{array}{l}
m_1(r),\quad r\in(0,z_{j_1}],\\
-m_1(r),\quad r\in(z_{j_1},a).
\end{array}\right.
$$
Then $m_2(r)\in C^1([0,1])$ with
 $|m_1(r)-m_2(r)|\le 2\tau_1$. From Lemma \ref{lemma} we obtain
$$
\big|\lambda(s_1,m_2(r))-\lambda(s_1,m_1(r))\big|\le c^*(e^{\frac{8}{s_1}}-1)<\frac{1}{2},
$$
and
$$
\big|\lambda(s_1,m_2(r))-\lambda^{\mathcal{D}}\big|\le \big|\lambda(s_1,m_2(r))-\lambda(s_1,m_1(r))\big|+\big|\lambda(s_1,m_1(r))-\lambda^{\mathcal{D}}\big|\le 1.
$$
Moreover, $m_2(r)\in S_{\mathcal{N}}$ from Example \ref{EE}, and
then we obtain from Theorem \ref{NT} that
$$\lim_{s\to+\infty}\lambda(s,m_2(r))=\lambda^{\mathcal{N}}.$$
Thus, there exists $s_2>\max\{s_1,8\mbox{/}\ln(1+\frac{1}{3 c^*})\}$ such that $$\big|\lambda(s_2,m_2(r))-\lambda^{\mathcal{N}}\big|<\frac{1}{3}.$$
Let $\tau_2=s_2^{-2}$. By the continuity of $m_2(r)$, there exists $\delta(\tau_2)\in(0,\delta(\tau_1))$, such that $|m_2(r)|<\tau_2$, as $r\in(a-\delta(\tau_2),a)$.
 \par

We select $j_2>j_1$ such that $z_{j_2}>a-\delta(\tau_2)$, and define
$$
m_3(r)=\left\{
\begin{array}{l}
m_2(r),\quad r\in(0,z_{j_2}],\\
-m_2(r),\quad r\in(z_{j_2},a).
\end{array}\right.
$$
Thus, $m_3(r)\in C^1([0,1])$ with $|m_2(r)-m_3(r)|\le 2\tau_2$. From Lemma \ref{lemma}, we obtain
$$
\big|\lambda(s_2,m_3(r))-\lambda(s_2,m_2(r))\big|\le c^*(e^{\frac{8}{s_2}}-1)<\frac{1}{3}.
$$
Then, we have
$$
\big|\lambda(s_2,m_3(r))-\lambda^{\mathcal{N}}\big|\le \big|\lambda(s_2,m_3(r))-\lambda(s_2,m_2(r))\big|+\big|\lambda(s_2,m_2(r))-\lambda^{\mathcal{N}}\big|\le \frac{2}{3}.
$$
Moreover, we note that $m_3(r)=m_1(r)$ for $r\in(z_{j_2},a)$, and by the definition of $m_1(r)$, we have $m_3(r)\in S_{\mathcal{D}}$.
Thus, we obtain from Theorem \ref{DT} that
$$\lim_{s\to+\infty}\lambda(s,m_3(r))=\lambda^{\mathcal{D}}.$$
\par
 We can repeat the above procedure to obtain
 $\{m_{n}(r)\}_{n=1}^{+\infty}$ 
 satisfying
 \begin{equation}
\big|\lambda(s_{2n},m_{2n+1}(r))-\lambda^{\mathcal{N}}\big|<\frac{2}{2n+1}\label{E:limitNN},
\end{equation}
where $s_{2n}>\max\{s_{2n-1},8\mbox{/}\ln(1+\frac{1}{(2n+1)c^*})\}$ and
\begin{eqnarray}
\big|\lambda(s_{2n+1},m_{2n+2}(r))-\lambda^{\mathcal{D}}\big|\le \frac{1}{n+1},\label{nn1}
\end{eqnarray}
where $s_{2n+1}>\max\{s_{2n},8\mbox{/}\ln(1+\frac{1}{(2 n+2)c^*})\}$, respectively.
We emphasize that the sequences constructed above satisfy that $s_n\to+\infty$, $\tau_n\to 0$, $\delta(\tau_n)\to 0$ and $z_{j_n}\to a$ as $n\to+\infty$.

\textbf{Step 2.} Construction of $m^*(r)$.

We claim that the function sequence $\{m_n(r)\}_{n=1}^{+\infty}$ is Cauchy sequence in $C^1([0,1])$. Since $m_n(r)=m_n(1-r)$ and  $m_n(r)=0$ on $[a,b]$, we have
\begin{equation}\label{rx}
\|m_n(r)-m_k(r)\|_{C^1([0,1])}=\|m_n(r)-m_k(r)\|_{C^1([0,a])}.
\end{equation}
Without loss of generality, assume $k\ge n$. Then, from the construction of $\{m_n(r)\}_{n=1}^{+\infty}$,
we have
\begin{equation}\label{r1}
|m_n(r)|=|m_1(r)|,\quad |m_n'(r)|=|m_1'(r)|,
\end{equation}
and
\begin{equation}\label{r2}
m_n(r)=m_k(r),\quad \text{for } k\ge n, \ r\in [0,a-\delta(\tau_n)].
\end{equation}
Therefore, invoking \eqref{r1} and \eqref{r2}, we can estimate
\begin{align}\label{r3}
\|m_n(r)-m_k(r)\|_{C^1([0,a])}&=\|m_n(r)-m_k(r)\|_{C^1([a-\delta(\tau_n),a])}\nonumber\\
&\le\|m_n(r)\|_{C^1([a-\delta(\tau_n),a])}+\|m_k(r)\|_{C^1([a-\delta(\tau_n),a])}\nonumber\\
&=2\|m_1(r)\|_{C^1([a-\delta(\tau_n),a])}.
\end{align}
Since $m_1(r)= m_{\mathcal{D}}(r)$ (see Example \ref{EE}) and $\delta(\tau_n)\to 0$ as $n\to +\infty$, we obtain that
\begin{equation}\label{1rx}
\|m_n(r)-m_k(r)\|_{C^1([0,1])}\le2\|m_1(r)\|_{C^1([a-\delta(\tau_n),a])}\to 0,\quad n\to +\infty,
\end{equation}
 due to \eqref{rx} and \eqref{r3}. Thus, there exists $m^*(r)\in C^1([0,1])$ satisfying $m_n(r)\to m^*(r) $ in $C^1([0,1])$.

\textbf{Step 3.} The non-existence of limit of $\lambda(s,m^*(r))$ as $s\to+\infty$.

We will prove this by constructing two subsequences that converge to different limits.\\
By the definition of $\{m_n(r)\}_{n=1}^{+\infty}$ and $m^*(r)$, we obtain
$$
|m_{2n+2}(r)-m^*(r)|<2\tau_{2n+2}, n\in \mathbb{N}.
$$
Then, by  Lemma \ref{lemma}, we have
$$
|\lambda(s_{2n+1},m_{2n+2})-\lambda(s_{2n+1},m^*(r))|\le c^*\left(e^{\frac{8}{s_{2n+2}}}-1\right)\le \frac{1}{2n+3}.
$$
Combining with \eqref{nn1}, we obtain
$$
|\lambda(s_{2n+1},m^*(r))-\lambda^{\mathcal{D}}|\le \frac{1}{n+1}+\frac{1}{2n+3}<\frac{2}{n}.
$$
Thus, $\lim\limits_{n\to+\infty }\lambda(s_{2n+1},m^*(r))= \lambda^{\mathcal{D}} $. Similarly, from \eqref{E:limitNN}, we directly obtain that $\lim\limits_{n\to+\infty }\lambda(s_{2n},m^*(r))= \lambda^{\mathcal{N}} $, completing the proof.
\qed

As a result, through iterative folding  on the initial function \( m_1 \in \mathcal{S}_{\mathcal{D}}\), we construct a sequence \( \{m_n\}_{n=1}^{+\infty} \) whose limit forms the desired example \(m^*\). The resulting function
$m^*$ degenerates to zero on the interval $[a,b]$, while exhibiting decaying oscillations that fluctuate up and down around the horizontal axis, gradually diminishing to zero as $r$ approach
$a$ from the left and $b$ from the right. A graph of $m(r)$ is illustrated in Figure 1.
\begin{figure}
  \centering
  \includegraphics[width=10cm]{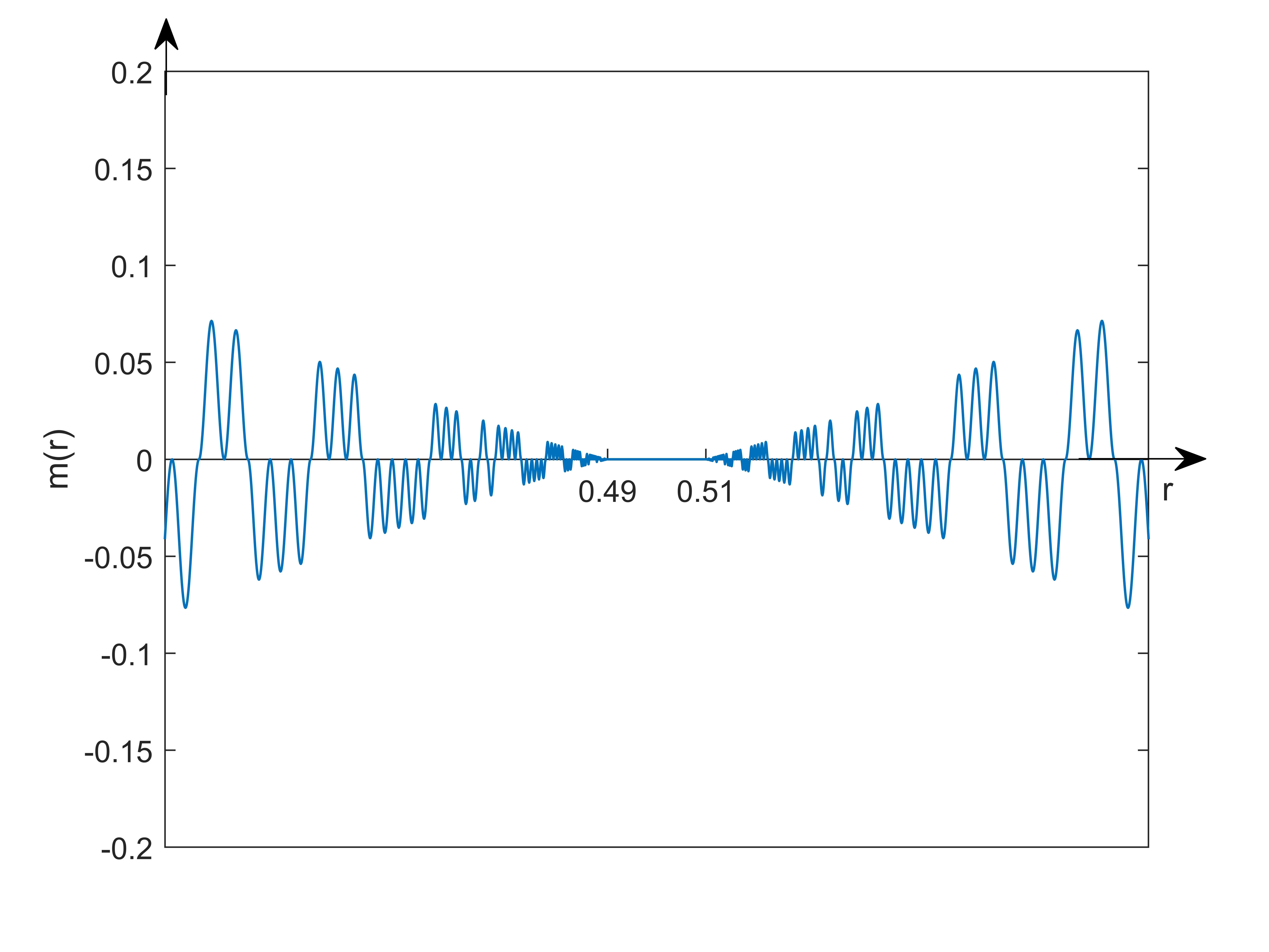}
  \vspace{-0.5cm}
\caption{An illustration of the graph of $m$ with $a=0.49$ and $b=0.51$.}\label{m}
\end{figure}


\section{An application}\label{S6}
We consider a one-dimensional single-species reaction-diffusion-advection model from \cite{book}.
\begin{equation}
\left\{
\begin{array}{llll}
u_t=u_{xx}+2sm_x(x)u_x+u(\sigma(x)-u),& 0<x<1,\ t>0,\\
u_x(0,t)=u_x(1,t)=0,& t>0,\\
u(x,0)=u_0(x)\ge,\not\equiv 0.
\end{array}\right.\label{E0}
\end{equation}
Here, $u(x,t)$ represents the population density. The function \( m(x) \in C^1([0,1])\) represents the driving force of external environmental factors on species advection, such as wind and water currents, and the parameter $s$ characterizes the advection strength. The space dependent function $\sigma\in C ([0,1])$ quantifies environmental suitability in space, integrating growth-promoting and resource-limiting factors. To compare the effects of different advection mechanisms, we fix the function $\sigma(x)$ and impose the following structural conditions for mathematical analysis.
\begin{assumption}\label{A:s}
For the constants $0<a<b<1$, $\sigma\in C ([0,1])$ satisfies the following conditions.

(1) The total mass of $\sigma(x)$ is positive: $\int_a^b \sigma(x)\mathrm{d} x>0$.

(2) The maximum $\sigma^{*}\triangleq\max_{x\in [a,b]} \sigma(x)$ satisfies $\frac{\pi^2}{(b-a)^2}-\sigma^{*}>0$.

(3) For small $\varepsilon>0$, $\sigma(x)$ is negative on $[a,a+\varepsilon)\cup(b-\varepsilon,b]$ and positive on $(a+\varepsilon,b-\varepsilon)$. Additionally, $\sigma(x)<\sigma(a)=\sigma(b)<-\frac{32}{3}(b-a-4\varepsilon)^{-2}$ for $x\in [0,a)\cup(b,1]$.
\end{assumption}
To proceed, we define $\tilde{\lambda}^{\mathcal{N}}$ and $\tilde{\lambda}^{\mathcal{D}}$ as the principal eigenvalues of the equations
$$
-\phi_{xx}-\sigma(x)\phi=\tilde{\lambda}\phi,\quad a<x<b,
$$
with Neumann boundary conditions ($\phi_x(a)=\phi_x(b)=0$) and Dirichlet boundary conditions ($\phi(a)=\phi(b)=0$), respectively. Then we have the following lemma.
\begin{lemma}\label{611}
For $\sigma$ satisfying Assumption \ref{A:s}, the following holds
$$\tilde{\lambda}^{\mathcal{N}}<0<\tilde{\lambda}^{\mathcal{D}}<-\sigma(x),\quad x\in[0,a]\cup[b,1].$$
\end{lemma}

\begin{proof}
Notice that the variational characterizations of $\tilde{\lambda}^{\mathcal{D}}$ and $\tilde{\lambda}^{\mathcal{N}}$ are same as the ones in \eqref{vd} and \eqref{vn} with $d=1$ and $c(x)\triangleq -\sigma(x)$, that is,
\begin{equation*}
\tilde{\lambda}^{\mathcal{D}}=\min_{\phi\in H_0^1((a,b))\atop\int_a^b\phi^2\mathrm{d} x=1}\int_a^b(|\phi'|^2-\sigma(x)|\phi|^2)\mathrm{d} x
\,\text{ and }\,\,
\tilde{\lambda}^{\mathcal{N}}=\min_{\phi\in H^1((a,b))\atop\int_a^b\phi^2\mathrm{d} x=1}\int_a^b(|\phi'|^2-\sigma(x)|\phi|^2)\mathrm{d} x.
\end{equation*}
Then testing with the constant function $\phi\equiv1$ on $[a,b]$ in the variational formula with condition (1) of Assumption \ref{A:s}, we obtain $\tilde{\lambda}^{\mathcal{N}}<0$.

Let $v$ be the eigenfunction corresponding to $\tilde{\lambda}^{\mathcal{D}}$. Since $v(a)=v(b)=0$, the Poincar\'{e} inequality  yields $\int_a^b |v'(x)|^2\mathrm{d} x \ge \frac{\pi^2}{(b-a)^2}\int_a^b v^2 \mathrm{d} x$. Substituting this into the variational formula for $\tilde{\lambda}^{\mathcal{D}}$ and using condition (2) in Assumption \ref{A:s}, we derive
$$
\tilde{\lambda}^{\mathcal{D}}\ge \frac{\int_a^b [\frac{\pi^2}{(b-a)^2}-\sigma(x)]v^2 \mathrm{d} x}{\int_a^b v^2 \mathrm{d} x}\ge \frac{\pi^2}{(b-a)^2}-\sigma^*>0.
$$

By condition (3) of Assumption \ref{A:s}, it suffices to show $\tilde{\lambda}^{\mathcal{D}}<-\sigma(a)$ to complete the proof. We construct a test function $\bar{v}$ as follows:
\begin{equation}\label{test1}
\bar{v}(x) = \left\{ \begin{array}{lcl}
0, && x\in[a,a+2\varepsilon),\\
\frac{x-a-2\varepsilon}{\delta}, && x\in[a+2\varepsilon,a+2\varepsilon+\delta),\\
1, & & x\in [a+2\varepsilon+\delta,b-2\varepsilon-\delta],\\
-\frac{x-b+2\varepsilon}{\delta}, & & x\in(b-2\varepsilon-\delta,b-2\varepsilon], \,\, \\
0,&& x\in[b-2\varepsilon,b], \end{array} \right.
\end{equation}
where $\delta=\frac{3}{8}(b-a-4\varepsilon)$ and $\epsilon>0$ is small enough. Substituting $\bar{v}$ into the variational formula for $\tilde{\lambda}^{\mathcal{D}}$ in \eqref{vd}, we get
\begin{align*}
\tilde{\lambda}^{\mathcal{D}}&\leq
\frac{\int_{a+2\varepsilon+\delta}^{a+2\epsilon} \delta^{-2} \mathrm{d} x-\int_{a+2\varepsilon}^{b-2\varepsilon} \sigma(x)\bar{v}^2\mathrm{d} x
+\int_{b-2\varepsilon-\delta}^{b-2\epsilon} \delta^{-2} \mathrm{d} x}{\int_{a+2\varepsilon}^{b-2\varepsilon} \bar{v}^2 \mathrm{d} x}\\
&=\frac{2\delta^{-1}-\int_{a+2\varepsilon}^{b-2\varepsilon} \sigma(x)\bar{v}^2\mathrm{d} x}{b-a-4\varepsilon-\frac{4}{3}\delta}\\
&\leq
\frac{2\delta^{-1}+(b-a-4\varepsilon-\frac{4}{3}\delta)\max_{x\in[a+2\varepsilon,b-2\varepsilon]}[-\sigma(x)]}
{b-a-4\varepsilon-\frac{4}{3}\delta}\\
&\leq 2\delta^{-1}(b-a-4\varepsilon-\frac{4}{3}\delta)^{-1}=\frac{32}{3}(b-a-4\varepsilon)^{-2}<-\sigma(a),
\end{align*}
where the last two inequalities follow from condition (3) in Assumption \ref{A:s}. This completes the proof.
\end{proof}

To see the impact of advection on population dynamics, we first consider a unidirectional advection term for comparison.

\textbf{Unidirectional advection ($m_x(x)\equiv-1$).}
In this case, the long-term population dynamics are determined by the stability of the trivial solution \( u\equiv0 \), which is determined by the principal eigenvalue of the associated  linearized equation of \eqref{E0}: 
\begin{equation}
\left\{
\begin{array}{l}
 -\phi_{xx}+2s\phi_x-\sigma(x)\phi=\lambda_1(s)\phi,\quad 0<x<1, \\
\phi_x(0 )=\phi_x(1 )=0,
\end{array}\right.\label{Ee}
\end{equation}
where $\lambda_1(s)$ is the principal eigenvalue with eigenfunction $\phi$.
The following theorem provides an explicit criterion for population persistence or extinction when the advection strength is large by [Theorem 1.1 \cite{CL}].
\begin{proposition}\label{p5}
For the eigenvalue problem \eqref{Ee} and model \eqref{E0} with $m_x(x)\equiv-1$, the following results hold. \\
(i) If $\sigma(0)>0$, then there exists $S>0$ such that for $s>S$, we have $\lambda_1(s)<0$, and \eqref{E0}  has a unique positive steady state $u^*$, which is globally asymptotically stable.\\
(ii) If $\sigma(0)<0$, then there exists $S>0$ such that for $s>S$, we have $\lambda_1(s)>0$, and \eqref{E0} has no positive steady state. All nonnegative solutions to \eqref{E0} decay exponentially to zero as $t\to \infty$.
\end{proposition}
 This proposition provides a criterion for the long-term dynamics of the equation with large advection. Simply speaking, when the advection rate \( s \) is sufficiently large, the species can persist if \( \sigma(0) > 0 \),
while will go extinct if \( \sigma(0) < 0 \). If we consider the interval as a river segment where $x=0$ corresponds to the upstream source and $x=1$ corresponds to the end of downstream.
Then Proposition \ref{p5} asserts that if the advection rate $s$ is large, the species can persist only if the some resource is supplied at the upstream source; otherwise, the species will go extinct in the river.
\begin{proof}
From Theorem 1.1 in \cite{CL}, we directly obtain that
$$
\lim_{s\to\infty}\lambda_1(s)=-\sigma(0).
$$
If $\sigma(0)>0$, there exists $S>0$, such that $\lambda_1(s)<0$, for $s>S$. According to Proposition 3.2 of \cite{book}, we can immediately derive (i). Similarly, if $\sigma(0)<0$, there exists $S>0$, such that $\lambda_1(s)>0$, for $s>S$, leading to (ii) in view of Proposition 3.1 of \cite{book}.
\end{proof}

\textbf{Case of complex advection (infinitely oscillating $m$).}
If the advection velocity field has complex structure like the potential function $m^*$ constructed in the proof of Theorem \ref{robin} in Section 5, the population dynamics can become complex. 

For the constructed potential function $m^*$, the eigenvalue problem of \eqref{E0} can be written as
\begin{equation}
\left\{
\begin{array}{l}
-\phi_{xx}-2sm^*_x \phi_x+c(x)\phi=\lambda_1(s;m^*)\phi(x),\quad 0<x<1,\\
\phi_x(0)=\phi_x(1)=0,
\end{array}\right.\label{E:example}
\end{equation}
where $c(x):= -\sigma(x)$ and $\lambda_1(s;m^*)$ is the principal eigenvalue of problem \eqref{E:example}.

We use Lemma \ref{611} and Theorem \ref{robin} to derive the following results that characterizing the population dynamics of \eqref{E0}, when the potential function $m^*$ has infinitely many oscillation.
\begin{proposition}\label{6}
For the eigenvalue problem \eqref{E:example} and model \eqref{E0} with $m(x)=m^*(x)$, under Assumption \ref{A:s}, the following results hold.  \\
(i) There exists a subsequence $\{s_i\}_{i=1}^{\infty}$ of $s\to\infty $, and an integer $I_0\in\mathbb{N}$, such that
 $\lambda_1(s_i;m^*)<0$ for  $i>I_0$. Along this sequence $\{s_i\}_{i>I_0}$, equation \eqref{E0}  has a unique positive steady state $u^*$, which is global asymptotically stable.\\
 (ii) There exists a subsequence $\{s_j\}_{j=1}^{\infty}$  of $s\to\infty $, and an integer $J_0\in\mathbb{N}$, such that
 $\lambda_1(s_j;m^*)>0$ for  $j>J_0$.
 Along this sequence $\{s_j\}_{j>J_0}$,
equation \eqref{E0} has no positive steady state, and all nonnegative solutions decay exponentially to zero as $t\to \infty$.
\end{proposition}
\begin{proof}


Applying Theorem \ref{robin} and  Lemma \ref{611}, we can find two subsequences $\{s_i\}_{i=1}^\infty$ and $\{s_j\}_{j=1}^\infty$ such that
$$
 \lim_{i\to\infty}\lambda_1(s_i;m^*)=\tilde{\lambda}^{\mathcal{N}}<0,\quad\quad  \lim_{j\to\infty}\lambda_1(s_j;m^*)=\tilde{\lambda}^{\mathcal{D}}>0.
$$
Therefore, there exist $I_0$ and $J_0$ such that
$\lambda_1(s_i;m^*)<0$ for  $i>I_0$, and $\lambda_1(s_j;m^*)>0$ for  $j>J_0$.
Applying Proposition 3.1 and Proposition 3.2 in \cite{book} along these subsequences $\{s_i\}_{i=1}^\infty$ and $\{s_j\}_{j=1}^\infty$, we could complete the proof.
\end{proof}
Proposition \ref{6} indicates that when the potential function is complex and the advection strength is strong, both persistence and extinction of the population are possible because the associated principal eigenvalue changes sign, which leads to a shift of the stability$\setminus$instability as $s$ increases.

For more clearly explanation, we consider the following example of $\sigma$.
\begin{example}
Let $a=\frac{1}{4}$, $b=\frac{3}{4}$. We construct a symmetric function $\sigma(x)$ on $[0,1]$, i.e.,  $\sigma(\frac{1}{2}-x)=\sigma(\frac{1}{2}+x)$, with the following piecewise form
\begin{equation*}
\sigma(x) = \left\{ \begin{array}{lcl}
-96, && x\in[0,\frac{1}{4}),\\
3072 x-864, && x\in[\frac{1}{4},\frac{9}{32}),\\
-512(x-\frac{1}{2})^2+\frac{49}{2}, & & x\in [\frac{9}{32},\frac{1}{2}].
\end{array} \right.
\end{equation*}
\end{example}
It is easy to verify that the above constructed $\sigma(x)$ satisfies Assumption \ref{A:s}. For this specific $\sigma(x)$, when $m_x=-1$, Proposition \ref{p5} shows that the population extinction occurs if the advection strength is large enough; while when $m=m^*$, Proposition \ref{6} asserts that the population dynamic become indeterminable for large advection strength, where either persistence or extinction is possible.

\bibliographystyle{abbrv}
\bibliography{bib}
\end{document}